\documentclass [11pt]{article}
 \usepackage{amssymb}
 \usepackage{amsmath,amsthm}
 \usepackage{enumitem}
 \usepackage{mathrsfs}
 \usepackage{color}
\usepackage[numbers,sort&compress]{natbib}
 \usepackage{graphicx}
\usepackage{tikz}
 \usepackage[colorlinks=true]{hyperref}

\normalsize

\renewcommand{\dfrac}[2]{\lower0.15ex\hbox{\large$\textstyle\frac{#1}{#2}$}}

\newtheorem{theorem}{Theorem}[section]
\newtheorem{lemma}[theorem]{Lemma}
\newtheorem{conjecture}[theorem]{Conjecture}

\newtheorem{corollary}[theorem]{Corollary}

\numberwithin{equation}{section}

\begin{document}
\title
{\bf Note on the random triadic process II  
\thanks{ The work was partially supported by NSFC (No.12071274).}}

\author
{{Fang Tian$^1$\qquad Yiting Yang$^2$}\\
{\small $^1$Department of Applied Mathematics}\\%
{\small Shanghai University of Finance and Economics, Shanghai 200433, China} \\
{\small\tt tianf@mail.shufe.edu.cn}\\[1ex]
{\small $^2$School of Mathematical Sciences}\\
{\small Tongji University, Shanghai 200092, China}\\
{\small\tt ytyang@tongji.edu.cn}
}

\date{}
 \maketitle

\begin{abstract}

For a fixed integer $r\geqslant 3$, the random $r$-generalized triadic process starts from
a minimum weakly $K_r$-saturated graph, namely the complete $(r-1)$-partite graph 
$K_{1^{r-2},n-r+2}$ with $r-2$ singleton parts
and one part of size $n-r+2$. At each step of the
process a new edge $\boldsymbol{e}$ is added with probability $p$, 
provided there exists a copy of $B_{r-2}$
in $K_n$ for which $\boldsymbol{e}$ is its spine and the only edge 
that has not yet been added to the current graph, where $B_{r-2}$
 denotes a triangular book graph consisting of $r-2$ triangles sharing a single
 common edge known as the spine of the book. 
 When $r=3$, this is the random triadic process. 
 Kor\'{a}ndi et al. proved a sharp threshold probability ${p}= \frac 12n^{ - \frac 12}$ for 
 the propagation of the random triadic process.
 Tian and Yang complemented those results by
 considering the final size of the random $r$-generalized 
triadic process when ${p}=o( n^{- \frac 12}\log^{ \alpha(3-r)} n)$
for any fixed constant $\alpha> \frac 12$. 
In this note, we further consider the sharp threshold 
probability $\tilde{p}_r$ for the propagation 
of the random $r$-generalized triadic process
 when $r\geqslant 4$. We prove that
$\frac 12 n^{- \frac 12}<
 \tilde{p}_r\leqslant c\, n^{- \frac 12}\log^{ \frac 12} n$
for some positive constant  $c$, and conjecture that 
$\tilde{p}_r= c_r\, n^{ - \frac 12}$, where 
$c_r= \frac 12(r-2)!^{ \frac{1}{2r-4}}\bigl( \frac{2r-5}{r-2}\bigr)^{ \frac{2r-5}{2r-4}}$.
\end{abstract}



\section{Introduction}

Bollob\'{a}s~\cite{bollo68} introduced the concept of \textit{weak saturation}, 
which gives rise to a process called $H$-\textit{bootstrap percolation}, starting from a graph $G$. 
Graphs are defined with respect to their edge sets. Given
a fixed template graph $H$ and an initial spanning subgraph $G\subseteq K_n$, 
the process evolves as follows.
 At each step of the process a new edge $\boldsymbol{e}$ 
 is added, provided there exists a copy of $H$ in $K_n$ for which
$\boldsymbol{e}$ is the only edge that has not yet been added to $G$,  where $K_n$
is the complete graph on the same set of vertices.
If no such edge exists, the process terminates. 
The closure of $G$ under the $H$-bootstrap process
is denoted by $\langle G\rangle_H$. When $\langle G\rangle_H=K_n$,
the subgraph $G$ is said to be \textit{weakly $H$-saturated}
in the terminology of Bollob\'{a}s~\cite{bollo68},
or equivalently, it \textit{$H$-percolates}. Bootstrap percolation 
is used to model evolving networks in various
real-world phenomena of interest. Although the 
dynamics are local, they can lead to global behavior. 
The minimum number of edges in a weakly $H$-saturated subgraph is 
called the \textit{weak $H$-saturation number} and 
is denoted by wsat$(K_n,H)$. For simplicity,
 wsat$(K_n,H)$ is written as wsat$(n,H)$. 
Lov\'{a}sz~\cite{lovasz77} proved that wsat$(n,K_r)=(r-2)n- \binom{r-1}{2}$
for $n\geqslant r\geqslant 2$, settling a conjecture of Bollob\'{a}s~\cite{bollo68}.
The result was also proved by Frankl~\cite{frankl82}, Kalai~\cite{kalai84},
Alon~\cite{alon1985} and Yu~\cite{yu1993}.

Balogh, Bollob\'{a}s and Morris~\cite{baloghbollo11} were the first to analyse the behaviour
of the $H$-bootstrap process relative to random initial edges, that is, when
the initial graph $G$ is chosen to be $\mathcal{G}({n,p})$,
 where $\mathcal{G}({n,p})$ is the Erd\H{o}s-R\'{e}nyi random graph obtained by choosing
each edge of $K_n$ independently with probability $p$. 
They defined the critical \textit{threshold} for $H$-bootstrap 
percolation on $K_n$ to be the point %
\begin{align*}
p_c(n,H)&:=\inf \Bigl\{p\in [0,1]: \mathbb{P}\bigl(\langle \mathcal{G}({n,p})\rangle_H=K_n\bigr)\geqslant \frac12\Bigr\},
\end{align*}
 at which $\mathcal{G}({n,p})$ is likely to $H$-percolate.
Bollob\'{a}s and Thomason~\cite{bollobas87} proved every non-trivial
monotone graph property has a threshold.
This is a threshold function of $H$-bootstrap percolation.

Balogh et al.~\cite{baloghbollo11}
initiated the study of $p_c(n,H)$ primarily focusing on the case $H=K_r$.
Note that all graphs $K_2$-percolate so trivially $p_c(n,K_2)=0$.
The classical result of Erd\H{o}s and R\'{e}nyi implies that $p_c(n,K_3)=n^{-1}\log n(1+o(1))$,
where $o(1)$ denotes a quantity that tends to zero
as $n\rightarrow \infty$. The threshold of $p_c(n,K_4)$  was 
identified up to a constant factor in~\cite{baloghbollo11},
and the existence of a sharp threshold at $p_c(n,K_4)= (3n\log n)^{-1/2}(1+o(1))$
follows from the works in~\cite{angel2018,angel2021,koles2022}. 
In a very recent breakthrough, Bartha, Kolesnik, Kronenberg
and Peled~\cite{bartha2025} established $p_c(n,K_r)= (\gamma n)^{-1/\lambda}(1+o(1))$ 
for every $r\geqslant 5$, where
$v_H$ and $e_H$ denote the numbers of vertices and 
edges for a graph $H$, and $\lambda= { (e_H-2)}/({v_H-2})$. 
The characterizations of $p_c(n,H)$ for some other types of template graphs $H$
 have been studied in~\cite{baloghbollo11}
and elsewhere, e.g. in~\cite{bartha20241,bartha20242,bayra2022,bidgoli21}.

Kor\'{a}ndi, Peled and Sudakov~\cite{korandi16} studied a related
process, that is, the \textit{random triadic process}.
The process starts with a minimum weakly $K_3$-saturated graph, 
a star $K_{1,n-1}$ on the vertex set $[n]$,
and repeatedly adds an edge $\boldsymbol{e}=\{x,y\}$ with probability $p$ provided 
that the end-vertices $x$ and $y$ already have a common neighbour in the graph. 
All existing edges are preserved in the graph.
The process is said to \textit{propagate} if all the edges are added to the graph eventually.
This process can be regarded as a class of relaxed dual problem of
the random bootstrap $K_3$-percolation process in~\cite{baloghbollo11}.
Their duality arises from an inversion of where randomness is introduced
and how edge expansion proceeds. At each step, all eligible edges are subject to
independent trials with probability $p$. The state distribution of 
the graph only depends on the number of iterative steps instead of the internal
checking order. Hence, the order of added edges 
does not affect the evolution of the graph process.

Let $\mathcal{P}$ be a monotone increasing graph property.
For a discrete random greedy graph process,
at each step, all eligible edges complying with some rule are 
independently added with probability $p$. A function $\tilde{p}=\tilde{p}(n)$ 
is a sharp threshold for $\mathcal{P}$ in the process,
for every fixed $\epsilon>0$, if the probability that the 
process satisfies  $\mathcal{P}$ converges to
 one when $p/\tilde{p}\geqslant 1+\epsilon$, while
 it converges to zero when $p/\tilde{p}\leqslant 1-\epsilon$. 
 The graph property of propagation is increasing 
and therefore has a threshold, as shown by Bollob\'{a}s and Thomason~\cite{bollobas87}. 
Kor\'{a}ndi et al.~\cite{korandi16} proved a sharp threshold probability for 
 the propagation of the random triadic process.
\begin{theorem}[Theorem~1.1 in~\cite{korandi16}]
If $p= c n^{ - \frac 12}$ for some positive constant $c$, 
then with high probability, the triadic process propagates 
when $c> \frac 12$ and stops at $O(n^{\frac 32})$ edges when $c< \frac 12$.
\end{theorem} 
\noindent As usual, an event is said to occur  with high probability (\textit{w.h.p.} for short), 
if the probability that it holds tends to 1 when $n\rightarrow\infty$.

 Recently,  for an approximately $d$-regular $\hbar$-uniform hypergraph $\mathcal{H}$
satisfying some mild degree and codegree conditions,
 Morrison and Noel~\cite{morr2021}{\color{red},}
motivated by Balogh et al.~\cite{baloghbollo11} and Kor\'{a}ndi
et al.~\cite{korandi16}{\color{red},} were the first to establish the sharp 
threshold result for bootstrap percolation 
in a random hypergraph using differential equation methods. 
Their model starts with random initial infected
vertices, where each vertex of 
$\mathcal{H}$ is included independently with probability $q$. 
At each step, a healthy vertex becomes infected independently with probability 
$p$ if there exists a hyperedge of $\mathcal{H}$ in which it
is the unique healthy vertex. Their analysis is based on the sufficient condition that both 
$p$ and $q$  are of the same order $d^{- 1/(\hbar-1)}$.
Given graphs $G$ and $F$, let $\mathcal{H}_{G,F}$ be a hypergraph in which each vertex corresponds to an edge
of $G$ and the hyperedges are precisely the edge sets of copies of $F$ in $G$.
Applying the result in~\cite{morr2021} to 
the special case  $\mathcal{H}=\mathcal{H}_{K_n,K_r}$,
it follows that  the expected number of initial vertices
needed for percolation is $n^{ \frac{2r}{r+1}}$ when 
$p$ is of the order of $n^{- \frac{2}{r+1}}$,
which far exceeds wsat$(n,K_r)=(r-2)n- \binom{r-1}{2}$. 
In fact,  the results in~\cite{morr2021} are valid for $\mathcal{H}_{K_n,F}$
provided $F$ is a strictly $2$-balanced graph. 
The authors stressed that alternative approaches 
 are required to address more general forms of 
$p$ and $q$.

For a fixed integer $r\geqslant 3$, 
the random $r$-generalized triadic process starts from
a minimum weakly $K_r$-saturated graph, a complete $(r-1)$-partite 
graph $K_{1^{r-2},n-r+2}$ on the vertex set $[n]$
with $(r-2)$ singleton parts and one part of size $n-r+2$. 
 At each step of the process a new edge $\boldsymbol{e}=\{x,y\}$ 
 is added with probability $p$, provided there exists a copy of $B_{r-2}$
in $K_n$ for which $\boldsymbol{e}$ is its spine and the only 
edge that has not yet been added to the current graph, where $B_{r-2}$
 denotes a triangular book graph consisting of $r-2$ triangles sharing a single
 common edge known as the spine of the book. In simple terms, this means 
 choosing two distinct vertices $x$ and $y$  uniformly at random if
they already have at least $r-2$ common neighbours in the graph
 and iteratively adding $\boldsymbol{e}=\{x,y\}$ as an edge with probability $p$.
Added edges are kept in the graph forever.
  The random triadic process is an abbreviation 
  for the random $3$-generalized triadic process. 
As with the random triadic process, the order 
in which edges are added has no impact on the evolution of the graph process.

For each fixed integer $r\geqslant 4$, in contrast with the case $r=3$, 
the process is closely related to the 
 well-known $r$-neighbour bootstrap percolation,
 where at each step, each non-infected vertex with at
 least $r$ infected neighborhoods becomes infected, and infected vertices remain infected forever.
As stated in~\cite{baloghbollo12},  this is similar
to $H$-bootstrap process with $H$ a star, except that a given copy of $H$ can only
be responsible for infecting its central vertex; it is thus a special case of
the natural `directed' version of $H$-bootstrap percolation.
The random $r$-generalized triadic process
when $r\geqslant 4$ is also a special case of the natural directed 
 version of $H$-bootstrap percolation with $H=B_{r-2}$, 
in which at each step, a given copy of $H$ 
 can only be responsible for adding its spine.
The  triangular book graph $B_{r-2}$ is 
 not a strictly $2$-balanced graph for $r\geqslant 4$.
The results of Morrison and Noel~\cite{morr2021}
cannot be applied to the random $r$-generalized triadic process for $r\geqslant 4$.

Tian and Yang~\cite{tiany24} considered the final size of the random $r$-generalized 
triadic process when $p=o( n^{- \frac 12}\log^{ \alpha(3-r)} n)$
for any fixed constant $\alpha> \frac 12$. Unlike the discussions in~\cite{korandi16,morr2021},
the use of differential equation methods is not 
required under this condition. We showed that 
the graph generated by the process essentially behaves like $\mathcal{G}(n,p)$,
and the number of added edges is dominated by the number of edges of $\mathcal{G}(n,p)$. 
\begin{theorem}[Theorem~1.1 in~\cite{tiany24}]
Let $r\geqslant 3$ be a fixed integer and $p=o( n^{- \frac 12}\log^{\alpha(3-r)} n)$
for any fixed constant $\alpha> \frac 12$. The final number of added edges 
in the random $r$-generalized triadic process, \textit{w.h.p.}, 
equals $ \frac {1}{2}n^{2}p(1\pm o(1))$  when $p=\omega(n^{-2})$.
\end{theorem}
\noindent For $r=3$, Theorem~1.2  complements 
Theorem~1.1 of~\cite{korandi16} on the random triadic process
to the case of $p=o( n^{- \frac 12})$.
\begin{corollary}
Suppose that $p=o( n^{- \frac 12})$. The final number of added edges 
in the random triadic process, \textit{w.h.p.}, equals
$ \frac {1}{2}n^{2}p(1\pm o(1))$  when $p=\omega(n^{-2})$.
\end{corollary}

In this note, we further complement Theorem~1.1 of~\cite{korandi16}
by studying the sharp threshold probability for 
 the propagation of the random $r$-generalized triadic process 
when $r\geqslant 4$. Differential equation methods are 
inapplicable here when $\tilde{p}=cn^{- \frac 12}$
for some constant $c$. The reason is that all variables 
except the codegrees in~\cite{korandi16} display concentration for $r=3$,
whereas none of the monitored variables exhibit concentration behavior when 
$r\geqslant 4$. All variables involved in the process are related to the codegrees when $r\geqslant 4$,
and a similar line of analysis to that of~\cite{korandi16} does not apply in our setting.
We adopt the notation from~\cite{tiany24} and  employ iterative applications of standard
concentration inequalities to prove Theorem~1.4 below.

\begin{theorem}
For each fixed integer $r\geqslant 4$, there exists a
positive constant $c$ such that, \textit{w.h.p.},
the random $r$-generalized 
triadic process propagates when $p\geqslant c n^{- \frac 12}\log^{ \frac 12} n$.
\end{theorem}
\noindent The graph property of propagation for the 
random $r$-generalized triadic process is increasing and
therefore has a threshold. Let $\tilde{p}_r$ denote the sharp propagation threshold
corresponding to the random $r$-generalized triadic process. 
By Theorem~1.1 in~\cite{korandi16}, it follows that
$\tilde{p}{_3}= \frac{1}{2}n^{- \frac 12}$.
Combining Theorem~1.4 with Theorem~1.1, 
for each fixed integer $r\geqslant 4$, 
we have the upper and lower bounds for 
the sharp propagation threshold in the random $r$-generalized
triadic process. 
\begin{corollary}
For each fixed integer $r\geqslant 4$, 
there exists a positive constant $c$  such that
the random $r$-generalized triadic process has a propagation threshold 
$\tilde{p}_r$ satisfying
\begin{align*}
\frac 12 n^{- \frac 12}<  \tilde{p}_r\leqslant c\, n^{- \frac 12}\log^{ \frac 12} n.
\end{align*}
\end{corollary}
\noindent

The remainder of the note is structured as follows. 
Notation used throughout the note is presented in Section 2.
We prove our main result Theorem~1.4 in Section 3. The final section concludes the note
and proposes a conjecture that $\tilde{p}_r= c_r\, n^{ - \frac 12}$ for 
$c_r= \frac 12(r-2)!^{ \frac{1}{2r-4}}( \frac{2r-5}{r-2})^{ \frac{2r-5}{2r-4}}$.

\section{Notation}

In view of Theorems~1.1 and~1.2, it suffices to assume that $r$ is a fixed integer
satisfying $r\geqslant 4$ below. Let $\mathbb{G}(0)$ be the complete $(r-1)$-partite 
graph $K_{1^{r-2},n-r+2}$ on the vertex set $[n]$,
where the $r-2$ singleton parts are labeled as $v_1,\cdots,v_{r-2}$,
respectively. Since the order in which candidates for edges 
are considered does not matter in the random
$r$-generalized triadic process, we may analyze the process in rounds. 
Let $u,v\in [n]$ be any pair of vertices 
and $W\subseteq [n]\backslash\{u,v\}$ be a subset of size $r-2$.
For an integer $i\geqslant 1$, denote the graph
after the $i$-th round by $\mathbb{G}(i)$, that is,
   ${\mathbb{G}}(i)$ depends on the outcomes of the process up
  to this round. At any round $i$  in the process, graphs are defined with 
  respect to edge sets on the vertex set $[n]$.
In $\mathbb{G}(i)$, the $r$-set $\{u,v\}\cup W$ is 
called  an \textit{open} $(r-2,2)$-walk at  $u$ or  $v$,
denoted as $uWv$, if the subgraph it induces contains
 a complete bipartite graph $K_{2,r-2}$ with vertex partitions $\{u,v\}$
 and  $W$, where at least one of the edges in this $K_{2,r-2}$ was added 
 in the $i$-th round and the pair $\{u,v\}$ has not yet been 
 added as  an edge. In one round of the process, 
expose all open $(r-2,2)$-walks of the form $uWv$ simultaneously and 
consider adding each corresponding pair $\{u,v\}$
as an edge independently with probability $p$.
Further define an open $(r-2,3)$-walk from $u$ to 
  $v$, denoted as $uWw'v$, to be obtained by attaching 
 the edge $w'v$ to an open $(r-2,2)$-walk $uWw'$ at the vertex $w'$;
 and an open $(r-2,4)$-walk between  $u$ and  $v$,
 denoted as $uW_1w'W_2v$, to be obtained by sticking together two open $(r-2,2)$-walks
$uW_1w'$ and $w'W_2v$ at $w'$, where $W_1,W_2\subseteq [n]\backslash\{u,v\}$ are
subsets of size $r-2$ and $w'\in [n]\backslash \{W_1\cup W_2\cup\{u,v\}\}$.

For any pair of vertices $u,v\in [n]$, 
let ${D}_u(i)$ be the set of neighbors of $u$,
${F}_u(i)$ be the set of open $(r-2,2)$-walks at  $u$,
 ${X}_{uv}(i)$ be the set of  new common neighbors
between $u$ and  $v$ that generated in the $i$-th round,
${Y}_{{uv}}(i)$ and 
${Z}_{uv}(i)$ be the sets of open $(r-2,3)$-walks
from $u$ to  $v$, and open $(r-2,4)$-walks between  $u$ and  $v$
in $\mathbb{G}(i)\backslash\{v_1,\cdots,v_{r-2}\}$, respectively.

For an event $A$ and a random variable $Z$ in an arbitrary 
probability space $(\Omega,\mathcal{F},\mathbb{P})$,
 $\mathbb{P}[A]$ and $\mathbb{E}[Z]$ denote the probability of  $A$
 and the expectation of  $Z$. We say that $X$ is a binomially distributed random 
 variable with parameters $n$ and $p$, or simply $X\sim\textbf{Bin}[n,p]$. 
 We will adopt the standard Landau notation. 
For sequences $\{a_n\}_{n=1}^\infty ,\{b_n\}_{n=1}^\infty$ of 
real numbers, write $a_n = O(b_n)$ if there exists a constant 
$K > 0$ such that for all $n$, $|a_n| \leqslant K  |b_n|$. 
Write $a_n= o(b_n)$ if eventually $b_n > 0$ and 
$\lim_{n \to \infty} \frac{a_n}{b_n}$ exists and equals $0$. 
Therefore, $a_n = b_n(1+o(1))$ means that $\lim_{n\to\infty} \frac{a_n}{b_n} = 1$. 
Denote $a_n = \Omega(b_n)$ if eventually $a_n > 0$ and $b_n = O(a_n)$. 
Further, $a_n = \Theta(b_n)$ if both $a_n = O(b_n)$ and $a_n = \Omega(b_n)$.
As usual, $[n]=\{1,\ldots,n\}$ and all asymptotics in this note are with respect to $n\to\infty$. 
All logarithms are natural, and the floor and ceiling signs are omitted
whenever they are not crucial.

\section{Main Results}

For $r=3$, in order to show $p_3= \frac 12 n^{- \frac 12}$
is the sharp propagation threshold, the proof in~\cite{korandi16}
consists of two phases. In the first phase, the authors made one step
 at a time and tracked the behavior of certain variables by differential equation methods 
in at most $n^2\sqrt{\log n}$ steps. Unlike in the first phase, they analyzed
the second phase of the process  in rounds such that it evolves into the 
complete graph in $O(\log \log n)$ rounds.  A similar analysis to that of~\cite{korandi16} 
does not work for $r\geqslant 4$.
The reason is that, when applying the differential equation method to 
the random $r$-generalized triadic process,
all variables except the codegrees in~\cite{korandi16} 
display concentration for $r=3$,
whereas none of the monitored variables exhibit concentration behavior when 
$r\geqslant 4$.
We can only use standard concentration inequalities to 
show that the process propagates in approximately $\frac{\log \log n}{\log (2r-4)}$ rounds
when $p\geqslant c\,n^{- \frac 12}\log^{ \frac 12}n$ for some positive constant $c$.
It suffices to suppose $p=c\,n^{- \frac 12}\log^{ \frac 12}n$
because the propagation property is increasing.
We make no attempts to optimize the constants below. 
There are many choices of them that can be balanced to satisfy certain inequalities.
Owing to the structure of $\mathbb{G}(0)$, it suffices,
when analyzing the  propagation threshold, to consider only
the vertices in $[n]\backslash\{v_1,\cdots,v_{r-2}\}$.

\begin{lemma} {\rm{(Chernoff bound~\cite{cher52})}}
Let  $X\sim\textbf{Bin}[n,p]$. Then, 
$\mathbb{P}[|X-np|>t]<2\exp\Bigl[- \frac{t^2}{3np}\Bigr]$
for any $0< t\leqslant np$. 
\end{lemma}

\begin{lemma} Let $r\geqslant 4$
be a fixed integer and $p=c\,n^{- \frac 12}\log^{ \frac 12}n$ for some positive constant $c$.
With probability at least $1- \frac 2n$, for any pair of vertices $u,v\in [n]\backslash\{v_1,\cdots,v_{r-2}\}$,
it follows that
\begin{enumerate}[label=(\arabic*).]
\item $|{D}_u(1)|=c\,(n\log n)^{ \frac 12}(1+o(1))$,
\item $|{X}_{uv}(1)|\in (\frac{c^2}{2}\log n, \frac{3c^2}{2}\log n)$,
\item $|{F}_u(1)|\geqslant \frac{c^{2r-4}}{2^{r-3}(r-2)!}(n\log^{r-2} n)(1+o(1))$, 
\item $|{Y}_{{uv}}(1)|\geqslant  \frac{c^{2r-3}}{2^{r-2}(r-2)!}(n^{ \frac 12}\log^{r- \frac 32} n)(1+o(1))$,
\item $|{Z}_{uv}(1)|\geqslant  \frac{c^{4r-8}}{2^{2r-5}(r-2)!^2}(n\log^{2r- 4} n)(1+o(1))$.
\end{enumerate}
\end{lemma}

\begin{proof}[Proof of Lemma~3.2]\, 
For any pair of vertices $u,v\in [n]\backslash\{v_1,\cdots,v_{r-2}\}$, 
we have that $\mathbb{G}(1)-\{v_1,\cdots,v_{r-2}\}$  is distributed as
 the random graph $\mathcal{G}(n-r+2,p)$.

Recall that ${D}_u(1)$ is the set of neighbors of $u$,
${X}_{uv}(1)$ is the set of  new common neighbors
that generated in the first round
between $u$ and  $v$ in $\mathbb{G}(1)\backslash\{v_1,\cdots,v_{r-2}\}$, respectively.
For any vertex $u\in [n]\backslash\{v_1,\cdots,v_{r-2}\}$,
the random variable $|{D}_u(1)|$ satisfies
$|{D}_u(1)|\sim\textbf{Bin}[n-r+1,p]$ and $\mathbb{E}[|{D}_u(1)|]=(n-r+1)p$. 
For any fixed $\epsilon>0$, applying Lemma 3.1 with $t= \epsilon \mathbb{E}[|{D}_u(1)|]$
 and noting that  $p=c\,n^{- \frac 12}\log^{ \frac 12}n$, we obtain 
$\mathbb{P}[||{D}_u(1)|-\mathbb{E}[|{D}_u(1)|]|>t]
<2\exp[-\epsilon^2(n-r+1)p]
< \frac{1}{n^2}$.
Taking a union bound over all vertices,
 with probability at most $ \frac 1n$, $|{D}_u(1)|=c\,(n\log n)^{ \frac 12}(1+o(1))$ 
 does not hold for some vertex $u\in [n]\backslash\{v_1,\cdots,v_{r-2}\}$. 

Similarly, for any pair of vertices $u,v\in [n]\backslash\{v_1,\cdots,v_{r-2}\}$, 
the random variable $|{X}_{uv}(1)|$ satisfies 
$|{X}_{uv}(1)|\sim\textbf{Bin}[n-r,p^2]$
and $\mathbb{E}[|{X}_{uv}(1)|]=(n-r)p^2= (c^2\log n)(1+o(1))$.
Upon substituting $t= \frac 13 c^2\log n$
into Lemma~3.1 with an appropriately chosen constant $c$, we further derive
$\mathbb{P}[||{X}_{uv}(1)|-\mathbb{E}[|{X}_{uv}(1)|]|>t]< 2\exp[- \frac{1}{27}(c^2\log n)(1+o(1))]
< \frac{1}{n^{3}}$.
Taking a union bound over all pairs of vertices,
with probability at most $ \frac 1n$, $|{X}_{uv}(1)|\in ( \frac{c^2}{2}\log n, \frac{3c^2}{2}\log n)$ 
does not hold for some pair of  vertices $u,v\in [n]\backslash\{v_1,\cdots,v_{r-2}\}$. 

In the sequel, assume that these events 
hold for all distinct vertices $u,v\in [n]\backslash\{v_1,\cdots,v_{r-2}\}$.
In order to obtain an open $(r-2,2)$-walk of the form  $uWv$
in ${F}_u(1)$, we first choose a neighbor of $u$, 
denoted as $w$,  in  $|{D}_u(1)|$ ways;
 and  choose a neighbor of  $w$ distinct from  $u$, denoted as $v$, 
 in $|{D}_w(1)|-1$ ways, from which we subtract the common 
 degrees of  $u$ and $w$; then we choose $(r-3)$ remaining common neighbors 
 of  $u$ and  $v$ from ${X}_{uv}(1)\backslash \{w\}$ in  $\binom{|{X}_{uv}(1)|-1}{r-3}$ ways. 
This procedure overcounts each set of $r-2$ common neighbors of  $u$ and $v$ exactly  $r-2$
 times. Hence, for any  $u\in [n]\backslash\{v_1,\cdots,v_{r-2}\}$, 
it follows that 
\begin{align*}
|{F}_u(1)|& \geqslant \frac{1}{r-2}\cdot |{D}_u(1)|\cdot \bigl(|{D}_w(1)|-O(|{X}_{uw}(1)|)\bigr)\cdot \binom{|{X}_{uv}(1)|-1}{r-3}\\
&\geqslant  \frac{c^{2r-4}}{2^{r-3}(r-2)!}(n\log^{r-2} n)(1+o(1)).
\end{align*}
In order to generate an open $(r-2,3)$-walk of the form $uWw'v$,
we first choose a neighbor of  $v$, denoted as $w'$, in $|{D}_v(1)|$
ways. Note that for any such choice of $w'$, 
in order to generate the open $(r-2,2)$-walk of the
form $uWw'$, we exclude the cases in which  $w'$ is same
as  $u$, or $w'$ is a common neighbor of  $u$ and  $v$. 
Combining these two cases, we count these choices of $w'$
in $O(|X_{uv}(1)|)$ ways, which equals $O(\log n)$ because 
$|{X}_{uv}(1)|\in ( \frac{c^2}{2}\log n, \frac{3c^2}{2}\log n)$.
We further take $(r-2)$ common neighbors of  $u$ and  $w'$
in $\binom{|{X}_{uw'}(1)|}{r-2}$ ways. It  follows that, 
for any pair of vertices $u,v\in [n]\backslash\{v_1,\cdots,v_{r-2}\}$,
\begin{align*}
|{Y}_{{uv}}(1)|&\geqslant \bigl(|{D}_v(1)|-O(\log n)\bigr)\cdot \binom{|{X}_{uw'}(1)|}{r-2}\\
&\geqslant  \frac{c^{2r-3}}{2^{r-2}(r-2)!}(n^{ \frac 12}\log^{r- \frac 32} n)(1+o(1)).
\end{align*}
In fact, we also have $|{Y}_{{uv}}(1)|\leqslant |{D}_v(1)|\cdot\binom{|{X}_{uw'}(1)|}{r-2}$
and $|{Y}_{{uv}}(1)|=\Theta(n^{ \frac 12}\log^{r- \frac 32} n)$.

 Finally, in order to obtain an open $(r-2,4)$-walk of the form  $uW_1w'W_2v$
in ${Z}_{uv}(1)$, we first choose an open $(r-2,2)$-walk at  $v$ in $|{F}_{v}(1)|$ ways,
denoted as $vW_2w'$ for some $w'\in [n]\backslash\{v_1,\cdots,v_{r-2}\}$. 
For each such choice $vW_2w'$,  in order to generate the open $(r-2,2)$-walk of the form $uW_1w'$,
we exclude the cases in which $w'$ is same as $u$, or $w'$ is a neighbor of  $u$.
Summing over these two bad cases, the total number of such 
 $w'$ is $O(|Y_{vu}(1)|)=\Theta(n^{ \frac 12}\log^{r- \frac 32} n)$.
Then we choose $(r-2)$ additional common neighbors 
 of $u$ and  $w'$ in  $\binom{|{X}_{uw'}(1)|}{r-2}$ ways. 
Hence, for any pair of vertices $u,v\in [n]\backslash\{v_1,\cdots,v_{r-2}\}$, 
\begin{align*}
|{Z}_{uv}(1)| &\geqslant \bigl(|{F}_v(1)|-O(|Y_{vu}(1)|)\bigr)\cdot \binom{|{X}_{uw'}(1)|}{r-2}\\
&\geqslant  \frac{c^{4r-8}}{2^{2r-5}(r-2)!^2}(n\log^{2r- 4} n)(1+o(1)).
\end{align*}
Combining the above, with probability at least $1- \frac 2n$, all these equations in Lemma~3.2 hold.
\end{proof}

For every pair of vertices $u,v\in [n]\backslash\{v_1,\cdots,v_{r-2}\}$,
define $\widetilde{X}_{uv}^{(r-2)}(1)$ to be
the set of open $(r-2,2)$-walks between $u$ 
and $v$ in $\mathbb{G}(1)\backslash\{v_1,\cdots,v_{r-2}\}$ 
if $u$ and $v$ are non-adjacent in $\mathbb{G}(1)$,
otherwise define $\widetilde{X}_{uv}^{(r-2)}(1)=\emptyset$.
For every non-adjacent pair $\{u,v\}$ in $\mathbb{G}(1)$,
based on the definition of $X_{uv}(1)$, 
we have $|\widetilde{X}_{uv}^{(r-2)}(1)|= \binom{|X_{uv}(1)|}{r-2}$.

\begin{lemma}With probability at least $1- \frac 2n$, 
all non-adjacent pairs $\{u,v\}$ in~$\mathbb{G}(1)$ have at least
$ \frac{n}{4^{2r-4}}$ vertices, denoted by $w$, such that 
\begin{align*}
|\widetilde{X}_{uw}^{(r-2)}(1)|\ \text{and }
 |\widetilde{X}_{vw}^{(r-2)}(1)|\geqslant  \frac{c^{2r-4}\log^{r-2}n}{6^{r-2}(r-2)!}(1+o(1)).
\end{align*}
\end{lemma}

\begin{proof}[Proof of Lemma~3.3]\, Assume that equations in Lemma~3.2
hold for all distinct vertices $u,v\in [n]\backslash\{v_1,\cdots,v_{r-2}\}$.
By the definition of ${Z}_{uv}(1)$, 
for every non-adjacent pair $\{u,v\}$ in~$\mathbb{G}(1)$,
 we have  
 \begin{align}
 |{Z}_{uv}(1)|&=\sum_{w\in [n]\backslash\{v_1,\cdots,v_{r-2},u,v\}}
 |\widetilde{X}_{uw}^{(r-2)}(1)|\cdot |\widetilde{X}_{vw}^{(r-2)}(1)|.
\end{align}
If $w$ is non-adjacent to both $u$ and $v$, then $|\widetilde{X}_{uw}^{(r-2)}(1)|= \binom{|{X}_{uw}(1)|}{r-2}$
and $|\widetilde{X}_{vw}^{(r-2)}(1)|= \binom{|{X}_{vw}(1)|}{r-2}$ in $\mathbb{G}(1)$.
It follows that 
\begin{align}
|\widetilde{X}_{uw}^{(r-2)}(1)|\ \text{and }|\widetilde{X}_{vw}^{(r-2)}(1)|\leqslant
\frac{3^{r-2} c^{2r-4}\log^{r- 2} n}{2^{r-2}(r-2)!}(1+o(1))
\end{align} by equation \textit{(2)} in
Lemma~3.2.

Suppose that there are strictly less than $ \frac{n}{4^{2r-4}}$ vertices $w$ such that 
\begin{align}
|\widetilde{X}_{uw}^{(r-2)}(1)|\cdot |\widetilde{X}_{vw}^{(r-2)}(1)|\geqslant
\frac{c^{4r-8}\log^{2r- 4} n}{2^{2r-4}\cdot (r-2)!^2}(1+o(1)).
\end{align}
Hence, by equations in~(3.1)-(3.3), we have 
\begin{align*}
|{Z}_{uv}(1)|&< \frac{c^{4r-8}n\log^{2r- 4} n}{2^{2r-4}\cdot (r-2)!^2}(1+o(1))+ 
 \frac{n}{4^{2r-4}}\cdot \frac{3^{2r-4} c^{4r-8}\log^{2r- 4} n}{2^{2r-4}(r-2)!^2}(1+o(1))\\
&<\frac{c^{4r-8}n\log^{2r- 4} n}{2^{2r-5}(r-2)!^2}(1+o(1)),
\end{align*}
contradicting equation ${(5)}$ in Lemma~3.2. Thus,
we have there are at least $ \frac{n}{4^{2r-4}}$ vertices $w$ satisfying (3.3).
Again by equation~(3.2), any vertex 
$w$ satisfying~(3.3) fulfills the requirements of Lemma 3.3.
\end{proof}

\begin{lemma} With probability at least $(1- \frac{2}{n})^2$, 
for each non-adjacent pair $\{u,v\}$ in~$\mathbb{G}(2)$, there exists a
 constant $0<\gamma<  \frac{1}{24^{2r-4}}$ such that
\begin{align*}
|X_{uv}(2)|\geqslant \frac{\gamma\cdot c^{4r-6}\log^{2r-3}n}{(r-2)!^2}(1+o(1)).
\end{align*}
\end{lemma}

\begin{proof}[Proof of Lemma~3.4] 
In the following, assume that the conditions in Lemma~3.3 hold for 
each  non-adjacent pair in~$\mathbb{G}(1)$.
Fix a non-adjacent pair $\{u,v\}$ in~$\mathbb{G}(1)$. 
In order to obtain the lower bound of $|X_{uv}(2)|$,
it suffices to consider these vertices $w$ such that the two inequalities
$|\widetilde{X}_{uw}^{(r-2)}(1)|$
and  $|\widetilde{X}_{vw}^{(r-2)}(1)|$ hold in 
Lemma~3.3. For each such $w$, it becomes a common neighbor of  $u$ and $v$ 
in $\mathbb{G}(2)$ with probability at least 
\begin{align*}
&\Bigl[1- \Bigl(1- cn^{- \frac 12}\log^{ \frac 12}n\Bigr)^{ \frac{c^{2r-4}\log^{r-2}n}{6^{r-2}\cdot(r-2)!}}\Bigr]^2
= \frac{c^{4r-6}n^{-1} \log^{2r-3}n}{6^{2r-4}\cdot(r-2)!^2}(1+o(1)).
\end{align*} 
Events associated with distinct vertices $w$ are independent. %
Hence, based on Lemma~3.3, $|X_{uv}(2)|$ 
stochastically dominates a random variable $X$ with
\begin{align*}
X\sim\textbf{Bin}\Bigl[ \frac{n}{4^{2r-4}}, \frac{c^{4r-6}n^{-1} \log^{2r-3}n}{6^{2r-4}\cdot(r-2)!^2}(1+o(1))\Bigr]
\end{align*}
and 
\begin{align*}
\mathbb{E}[X]= \frac{c^{4r-6} \log^{2r-3}n}{24^{2r-4}\cdot(r-2)!^2}(1+o(1)).
\end{align*}
Applying Lemma~3.1 again, there exists a constant  $0<\gamma<  \frac{1}{24^{2r-4}}$ such that 
\begin{align*}
&\mathbb{P}\Bigl[|X_{uv}(2)|<
\frac{\gamma\cdot c^{4r-6}\log^{2r-3}n}{(r-2)!^2}(1+o(1))\Bigr]\\
&< \mathbb{P}\Bigl[X<
\frac{\gamma\cdot c^{4r-6}\log^{2r-3}n}{(r-2)!^2}(1+o(1))\Bigr]\\
&= \mathbb{P}\Bigl[X<
(1-24^{2r-4}\gamma)\mathbb{E}[X]\Bigr]\\
&=\exp[-\Omega (\log^{2r-3}n)]\\
&< \frac{1}{n^3}. %
\end{align*}
Taking a union bound over all pairs of vertices 
and using the law of total probability, we complete the proof of Lemma~3.4. 
\end{proof}

To give a lower bound on $|X_{uv}(i)|$ 
for each non-adjacent pair $\{u,v\}$ in~$\mathbb{G}(i)$,
define a sequence $x_i$ for $2\leqslant i\leqslant M+1$ as 
\begin{align}
x_i&=\begin{cases} \Bigl(\frac{\gamma\cdot c^{4r-6}}{(r-2)!^2}\log^{2r-3}n\Bigr)^{(2r-4)^{i-2}},&i=2,\ldots, M,\\
(n\log n)^{ \frac{1}{2r-4}},&i=M+1;\end{cases}
\end{align} 
where $M$ is defined as $M=\min\{i:\, x_i^{2r-4}\geqslant (n\log n)^{ \frac{1}{2r-4}}\}$.
This implies that $M= \frac{\log\log n}{\log(2r-4)}(1+o(1))$.
\begin{lemma}With probability at least $(1- \frac 2n)^{i}$, 
for each non-adjacent pair $\{u,v\}$ in~$\mathbb{G}(i)$
\begin{align*}
|X_{uv}(i)|&\geqslant x_i(1+o(1)),
\end{align*} 
where $2\leqslant i\leqslant M+1$.
\end{lemma}

\begin{proof}[Proof of Lemma~3.5]
We prove this by induction on the number of rounds. Obviously, the equation
holds when $i=2$ by Lemma~3.4. Assume  $2\leqslant i<M+1$. 
In the following, we condition on the event that $|X_{uv}(i)|\geqslant x_i(1+o(1))$ 
holds  for each non-adjacent pair $\{u,v\}$ in~$\mathbb{G}(i)$,
and proceed to bound the probability that the event
fails for some non-adjacent pair
in the $(i+1)$-round. 

Let $p_{i+1}=1- (1-p)^{\binom{x_i(1+o(1))}{r-2}}= \frac{px_i^{r-2}}{(r-2)!}(1+o(1))$, where
$p=cn^{- \frac 12}\log^{ \frac 12}n$. 
Since all these open $(r-2,2)$-walks in $\mathbb{G}(i)$ are exposed simultaneously, 
each non-adjacent pair has at least $\binom{x_i(1+o(1))}{r-2}$ independent 
opportunities to be added with probability $p$ 
in the $(i+1)$-th round, which implies that every non-adjacent pair is added
in $\mathbb{G}(i+1)$ independently with probability at least $p_{i+1}$. 
Since all existing edges are preserved in the graph,
there exists a coupling such that  $\mathcal{G}(n,p_{i+1})$ is a subgraph of $\mathbb{G}(i+1)$
almost surely.

Fix a  non-adjacent pair $\{u,v\}$ in $\mathbb{G}(i+1)$.
The number of new common neighbors of  $u$ and $v$,
that is $|X_{uv}(i+1)|$, stochastically dominates
a random variable $X$ with $X\sim\textbf{Bin}[ n-2, p_{i+1}^2]$. It is easy to derive
\begin{align*}
\mathbb{E}[X]&= \frac{np^2 \cdot x_{i}^{2r-4}}{(r-2)!^2}(1+o(1))\geqslant \frac{c^2\log n \cdot x_{i+1}}{(r-2)!^2}(1+o(1))
\end{align*}
 by equation in~(3.4) and $M=\min\{i:\, x_i^{2r-4}> (n\log n)^{ \frac{1}{2r-4}}\}$.
By Lemma~3.1, we again obtain 
\begin{align*}
&\mathbb{P}\bigl[|X_{uv}(i+1)|<
x_{i+1}(1+o(1))\bigr]\\
&<\mathbb{P}\bigl[X<
x_{i+1}(1+o(1))\bigr]\\
&=\exp[-\Omega(\log n\cdot x_{i+1})]\\
&< \frac{1}{n^3}.
\end{align*}
Applying a union bound over all vertex pairs together 
and using the law of total probability,
we complete the proof of Lemma~3.5. 
\end{proof}

\begin{proof}[Proof of Theorem~1.4]
We claim that $\mathbb{G}(M+2)$ is the complete graph. 
For any non-edge $\boldsymbol{e}$ in $\mathbb{G}(M+1)$, denoted as $\boldsymbol{e}=\{u,v\}$,
by Lemma~3.5 and the definition of $X_{uv}(M+1)$,
it follows that the probability that the edge $\boldsymbol{e}$
is absent from $\mathbb{G}(M+2)$ is bounded above by
\begin{align*}
&(1-p)^{\binom{x_{M+1}(1+o(1))}{r-2}}\\
&\leqslant \exp\Bigl[-\frac{c\log n}{(r-2)!}\Bigr](1+o(1))\\
&< \frac{1}{n^3},
\end{align*}
where the final inequality follows from $x_{M+1}=(n\log n)^{ 1/(2r-4)}$
in~(3.4) with a suitably chosen constant $c>0$.
Taking a union bound over all pairs of vertices and using the law of total probability, 
we obtain that $\mathbb{G}(M+2)$ is  complete 
with probability at least $(1- \frac 2n)^{M+2}\rightarrow 1$,
where $M=\Theta(\log\log n)$. 
\end{proof}

Let integers $r_1$ and $r_2$ satisfy $3\leqslant r_1<r_2$.
For each $i\in \{1,2\}$, let $(\Omega_{r_i}, \mathcal{F}_{r_i}, \mathbb{P}_{r_i})$
denote the probability space associated with the random $r_i$-generalized triadic
process. We build a product coupling space $(\Omega, \mathcal{F}, \mathbb{P})$
with $\Omega=\Omega_{r_1}\times \Omega_{r_2}$ the standard Cartesian product of two sample
spaces, $\mathcal{F}=\mathcal{F}_{r_1}\otimes\mathcal{F}_{r_2}$
the minimal $\sigma$-algebra generated by all $A\times B$ with $A\in \mathcal{F}_{r_1}$
and $B\in \mathcal{F}_{r_2}$; a coupling probability measure $\mathbb{P}$
on $(\Omega, \mathcal{F})$ that preserves the original 
distribution of each individual process, 
that is, $\mathbb{P}(A\times\Omega_{r_2})=\mathbb{P}_{r_1}(A)$
and $\mathbb{P}(\Omega_{r_1}\times B)=\mathbb{P}_{r_2}(B)$
for all $A\in \mathcal{F}_{r_1}$ and $B\in \mathcal{F}_{r_2}$.

Let $A_{r_1}$ and $B_{r_2}$ be the events that the random $r_1$-generalized triadic process
propagates and the random $r_2$-generalized triadic process propagates, respectively.
Define coupled versions of these events as $\tilde{A}_{r_1}={A}_{r_1}\times \Omega_{r_2}$
and $\tilde{B}_{r_2}= \Omega_{r_1}\times{B}_{r_2}$
such that $\tilde{B}_{r_2}\subsetneq\tilde{A}_{r_1}$.
By monotonicity of probability measures $\mathbb{P}(\tilde{B}_{r_2})<
\mathbb{P}(\tilde{A}_{r_1})$, we have $\mathbb{P}_{r_2}({B}_{r_2})< \mathbb{P}_{r_1}({A}_{r_1})$.
We further deduce the propagation threshold $\tilde{p}_{r_2}> \tilde{p}_{r_1}$.
By Theorems~1.1 and~1.4, we obtain Corollary~1.5.

\section{Conclusion}

For the random bootstrap $K_r$-percolation process, 
the percolation threshold is approximately equal to
 $n^{-1}\log n $ for $r=3$, $(3n\log n)^{-1/2}$ for $r=4$ 
 in~\cite{angel2018,angel2021,koles2022}, 
and $(\gamma n)^{-1/\lambda}$ for some constant $\lambda$ and all $r\geqslant 5$
in~\cite{bartha2025}. The random $r$-generalized triadic process 
can be regarded as a class of  dual problems 
associated with random bootstrap $K_r$-percolation. 
Its sharp propagation threshold was shown to be 
$\frac 12 n^{ - \frac 12}$ for $r=3$ in~\cite{korandi16}.
As one special case in~\cite{morr2021}, it follows that the expected
number of initial vertices needed for percolation is 
of order $n^{ \frac{3}{2}}$ for $r=3$ when $p$ is 
of the order of $n^{- \frac{1}{2}}$, which far exceeds wsat$(n,K_3)=n- 1$. 
It is interesting to consider the sharp propagation threshold 
of the random $r$-generalized triadic process for $r\geqslant 4$.
In the present work, we establish that there exist constant
 $c$ such that its propagation threshold lies 
between $\frac 12 n^{- \frac 12}$ and $c\, n^{- \frac 12}\log^{ \frac 12} n$
for all $r\geqslant 4$.

Suppose $\tilde{p}_r=c\ n^{ -\frac 12}$ for some constant 
$c$. A pseudo-random heuristic, following
 the same analysis as~\cite{korandi16,morr2021},
predicts the evolution of variables 
and suggests that the expected trajectory of $|F_u(i)|$ 
for any vertex $u\in [n]\backslash \{v_1,\ldots,v_{r-2}\}$
is $f(t)\cdot n$, where $t= \frac{i}{n^2}$ and
 \begin{align*}
 f(t)&= \frac{(2c)^{2r-4}}{(r-2)(r-2)!}t^{2r-4}-2t+1.
 \end{align*}
  If $c> \frac 12(r-2)!^{ \frac{1}{2r-4}}( \frac{2r-5}{r-2})^{ \frac{2r-5}{2r-4}}$,
  then $f(t)$ has no positive root and the number of 
  open $(r-2,2)$-walks will be large enough that we 
  expect that the propagation occurs with high probability. 
  If $c< \frac 12(r-2)!^{ \frac{1}{2r-4}}( \frac{2r-5}{r-2})^{ \frac{2r-5}{2r-4}}$,
  then we expect the number of open $(r-2,2)$-walks to 
  become very small, suggesting that the propagation does not occur with high probability.
\begin{conjecture}
For each fixed integer $r\geqslant 4$, the sharp propagation threshold 
$\tilde{p}_r$ corresponding to the random 
$r$-generalized triadic process takes the form
$\tilde{p}_r=  c_r n^{ - \frac 12}$ for a constant $c_r$, where
\begin{align*}
c_r= \frac 12(r-2)!^{ \frac{1}{2r-4}}\Bigl( \frac{2r-5}{r-2}\Bigr)^{ \frac{2r-5}{2r-4}}.
\end{align*}
\end{conjecture}

\end{document}